\documentclass[a4paper,11pt,english]{article}
\usepackage[utf8]{inputenc}
\usepackage{babel,amsmath,amssymb,amsthm,graphicx,titlesec}
\usepackage[square,numbers,sort&compress]{natbib}
\usepackage[margin=2.5cm]{geometry}
\usepackage[font=small,labelfont=bf]{caption}

\titleformat{\section}{\normalfont\bfseries}{\thesection.}{1em}{}
\titlespacing*{\section}{0pt}{7ex plus 1ex minus .2ex}{2.5ex plus .2ex}

\let\P\undef
\let\Re\undef
\DeclareMathOperator{\P}{P}
\DeclareMathOperator{\Re}{Re}
\DeclareMathOperator{\E}{E}
\DeclareMathOperator{\I}{I}
\DeclareMathOperator{\cov}{cov}
\DeclareMathOperator{\sgn}{sgn}
\DeclareMathOperator{\Law}{Law}
\newcommand{\FF}{{{}_2F_1}}
\newcommand{\F}{\mathcal{F}}
\newcommand{\C}{\mathrm{C}}
\newcommand{\R}{\mathbb{R}}
\renewcommand{\phi}{\varphi}
\renewcommand{\epsilon}{\varepsilon}
\renewcommand{\tilde}{\widetilde}

\newtheorem*{theorem}{Main theorem}

\title{A Bayesian sequential test\\
for the drift of a fractional Brownian motion}

\author{Alexey Muravlev\thanks{Steklov Mathematical Institute of Russian
Academy of Sciences, 8 Gubkina st., Moscow 119991, Russia. Email:
almurav@mi.ras.ru.}
\and
Mikhail Zhitlukhin\thanks{Steklov Mathematical Institute of Russian Academy
of Sciences, 8 Gubkina st., Moscow 119991, Russia. Email:
mikhailzh@mi.ras.ru. Corresponding author.} }

\begin{document}
\maketitle

\begin{abstract}
We consider a fractional Brownian motion with unknown linear drift such that
the drift coefficient has a prior normal distribution and construct a
sequential test for the hypothesis that the drift is positive versus the
alternative that it is negative. We show that the problem of constructing
the test reduces to an optimal stopping problem for a standard Brownian
motion, obtained by a transformation of the fractional one. The solution is
described as the first exit time from some set, whose boundaries are shown
to satisfy a certain integral equation, which is solved numerically.

\smallskip
\textit{Keywords:} sequential test, Chernoff's test, fractional Brownian
motion, optimal stopping.

\smallskip
\textit{AMS Subject Classification:} 62L10, 62L15, 60G40.
\end{abstract}

\section{Introduction}
Suppose one observes a fractional Brownian motion process (fBm) with linear
drift and unknown drift coefficient. We are interested in sequentially
testing the hypotheses that the drift coefficient is positive or negative.
By a sequential test we call a procedure which continues to observe the
process until a certain time (which generally depends on a path of the
process, so it is a stopping time), and then decides which of the two
hypotheses should be accepted. We consider a Bayesian setting where the
drift coefficient has a prior normal distribution, and we use an optimality
criterion of a test which consists of a linear penalty for the duration of
observation and a penalty for a wrong decision proportional to the true
value of the drift coefficient. The goal of  this paper is to describe the
structure of the exact optimal test in this problem, i.e. specify a
stopping time and a rule to choose between the two hypotheses.

The main novelty of our work compared to the large body of literature
related to sequential tests (for an overview of the field, see e.g.
\cite{L97,TNB}) is that we work with fBm. To the best of our knowledge, this is the first
non-asymptotic solution of a continuous-time sequential testing problem for this process.
Recall that the fBm is a Gaussian process, which generalizes the standard
Brownian motion (sBm) and allows for dependent increments; see the
definition in Section \ref{setup}. It was first introduced by Kolmogorov
\cite{K40} and gained much attention after the work of Mandelbrot and van
Ness \cite{MvN}. Recently, this process has been used in various models in
applied areas, including, for example, modeling of traffic in computer
networks and modeling of stock market prices and their volatility; a
comprehensive review can be found in the preface to the
monograph~\cite{M07}.

It is well-known that a fBm is not a Markov process, neither a
semimartingale except the case when it is a sBm. As a
consequence, many standard tools of stochastic calculus and stochastic
control (It\^o's formula, the HJB equation, etc.) cannot be applied in
models based on fBm. In particular, recall that a general method to
construct exact sequential tests, especially in Bayesian problems, consists
in reduction to optimal stopping problems for processes of sufficient
statistics; see e.g. Chapter~VI in the book~\cite{PS}. In the majority of
problems considered in the literature sufficient statistics are Markov
processes, and so the well-developed general theory of Markov optimal
stopping problems can be applied. On the other hand, optimal stopping
problems for non-Markov and non-semimartingale processes, like fBm, often
cannot be solved even by numerical methods.

Fortunately, in the problem we consider it turns out to be possible to
change the original problem for fBm so that it becomes tractable. One of the
key steps is a general transformation outlined in the note \cite{M13}, which
allows to reduce sequential testing problems for fBm to problems for
diffusion processes. It is achieved by integration of a certain kernel with
respect to an observable process and using the known fact that a fBm can
be expressed as an integral with respect to a sBm, and vice versa.

In the literature, the result which is most closely related to ours is the
sequential test proposed by H. Chernoff \cite{C1}, which has exactly the
same setting and uses the same optimality criterion, but considers only sBm.
For a prior normal distribution of the drift coefficient, Chernoff and
Breakwell \cite{C2,C3} found asymptotically optimal sequential tests when
the variance of the drift goes to zero or infinity. In the paper
\cite{ZM12}, we extended their result and
constructed an exact optimal test. An important step was a transformation of
the problem that reduced the optimal stopping problem for the sufficient
statistic process, as studied by Chernoff and Breakwell, to an optimal
stopping problem for a standard Brownian motion with nonlinear observation
cost. A similar transformation is used in the present paper as well, see
Section~\ref{reduction}.

Let us mention two other recent results in the sequential analysis of fBm,
related to estimation of its drift coefficient. Cetin, Novikov, and Shiryaev
\cite{CNS} considered a sequential estimation problem assuming a normal
prior distribution of the drift with a quadratic or a $\delta$-function penalty for a wrong
estimate and a linear penalty for observation time. They proved that in
their setting the optimal stopping time is non-random. Gapeev and Stoev
\cite{GS} studied sequential testing and changepoint detection problems for
Gaussian processes, including fBm. They showed how those problems can be
reduced to optimal stopping problems and found asymptotics of optimal
stopping boundaries. There are many more results related to fixed-sample
(i.e. non-sequential) statistical analysis of fBm. See, for example, Part~II
of the recent monograph \cite{T17}, which discusses statistical methods for
fBm in details.

Our paper is organized as follows. Section~\ref{setup} formulates the
problem. Section~\ref{reduction} describes a transformation of the original
problem to an optimal stopping problem for a sBm and introduces auxiliary
processes which are needed to construct the optimal sequential test. The
main result of the paper -- the theorem which describes the structure of the
optimal sequential test -- is presented in Section \ref{main}, together with
a numerical solution. Section~\ref{proofs} contains the proof. Some
technical details are in the appendix.

\section{Decision rules and their optimality}
\label{setup}
Suppose one observes the stochastic process
\[ 
Z_t = \theta t + B_t^H, \qquad t\ge 0, 
\]
where $B^H_t$, $t\ge0$, is a fractional Brownian motion (fBm) with known
Hurst parameter $H\in(0,1)$ and unknown drift coefficient $\theta$. Recall
that a fBm is a continuous zero-mean Gaussian process with the covariance function
\[ 
\cov(B_t^H, B_s^H) = \frac12 (s^{2H} + t^{2H} - |t-s|^{2H}), \qquad
t,s\ge0.
\] 
In the particular case $H=1/2$ this process is a standard Brownian motion
(sBm) and has independent increments; its increments are positively
correlated in the case $H>1/2$ and negatively correlated in the case
$H<1/2$. Except the case $H=1/2$, a fBm is not a Markov process, neither a
semimartingale.

We will consider a Bayesian setting and assume that $\theta$ is a random
variable defined on the same probability space as $B^H$, independent of it
and having a normal distribution with known mean $\mu\in\mathbb{R}$ and
known variance $\sigma^2>0$. 

It is assumed that neither the value of $\theta$, nor the value of $B_t^H$
can be observed directly, but the observer wants to determine whether the
value of $\theta$ is positive or negative based on the information conveyed
by the combined process $Z_t$. We will look for a sequential test for the
hypothesis $\theta>0$ versus the alternative $\theta\le 0$. By a sequential
test we call a pair $\delta=(\tau,d)$ which consists of a stopping time
$\tau$ of the filtration $\F_t^Z$, generated by $Z$, and an
$\F_\tau^Z$-measurable function $d$ assuming values $\pm1$. The stopping
time is the moment of time when observation is terminated and a decision
about the hypotheses is made; the value of $d$ shows which of them is
accepted.

We will use the criterion of optimality of a decision rule consisting in
minimizing the linear penalty for observation time and the penalty for a
wrong decision proportional to the absolute value of $\theta$. Namely, with
each decision rule $\delta$ we associate the risk
\begin{equation} 
R(\delta) = \E( \tau + |\theta| \I(d \neq \sgn(\theta))),
\label{risk}
\end{equation}
where $\sgn(\theta) = -1$ if $\theta\le 0$ and $\sgn(\theta)=1$ if $\theta>0$. 
The problem consists in finding $\delta^*$ that minimizes $R(\delta)$ over
all decision rules. Note that one can consider a more general setting when
the penalty for observation time is equal to $c\tau$ with some constant $c>0$ (or
the penalty for a wrong decision is $c|\theta|\I(d \neq \sgn(\theta))$), but
this case can be reduced to the one we consider by a change of the
parameters $\mu,\sigma$ (see \cite{C1}), and so we'll focus only on $c=1$.

This problem was proposed by H.~Chernoff in \cite{C1} for sBm, and we refer
the reader to that paper and the subsequent papers \cite{C2,C3,C4} for a rationale for this setting. Those
papers contain results about the asymptotics of the optimal
test and other its properties, including a comparison with Wald's sequential
probability ratio test. Our paper \cite{ZM12} contains a result which allows
to find an exact (non-asymptotic) optimal test by a relatively simple
numerical procedure.

\section{Reduction to an optimal stopping problem}
\label{reduction}
We will transform the problem of finding a decision rule minimizing
\eqref{risk} by eliminating the function $d$ from it and  reducing it to
an optimal stopping problem. On the first step, we'll obtain an optimal
stopping problem for a fBm. Then by changing time and space coordinates it
will be reduced to an optimal stopping problem for a sBm, which will allow
to apply well-developed methods to solve it.

From the relation
$|\theta| \I(d\neq\sgn(\theta)) = \theta^+\I(d=-1) + \theta^-\I(d=1)$, where
$\theta^+ = \max(\theta,0)$, $\theta^- = -\min(\theta,0)$, and the fact that
$d$ is $\F_\tau^Z$-measurable, one can see that the optimal decision rule
should be looked for among rules $(\tau,d)$ with
$d= \min(\E(\theta^- \mid \F_\tau^Z), \E(\theta^+ \mid \F_\tau^Z))$. Hence,
it will be enough to solve the optimal stopping problem which consists in
finding a stopping time $\tau^*$ such that
$R(\tau^*) = \inf_{\tau} R(\tau)$, where
\[
R(\tau) =  \E(\tau + \min (\E(\theta^- \mid \F_\tau^Z),
\E(\theta^+ \mid \F_\tau^Z)))
\]
(for brevity, we'll use the same notation $R$ for the functional associated
with a decision rule, and the functional associated with a stopping time).
Then the optimal decision rule will be $\delta^* = (\tau^*,d^*)$ with
$d^*= 1$ if $\E(\theta^- \mid \F_{\tau^*}^Z)< \E(\theta^+ \mid
\F_{\tau^*}^Z))$ and $d^*=-1$ otherwise.

Next we are going to transform the expression inside the expectation in
$R(\tau)$ to the value of some process constructed from a sBm. It is known
(see e.g. \cite{J06} and the earlier results \cite{NVV,MG}) that the following
process $B_t$, $t\ge0$, is a sBm, and the filtrations generated by $B_t$ and
$B_t^H$ coincide:
\begin{equation}
B_t = C_H \int_0^t K_H(t,s) d B_s^H
\label{B via BH}
\end{equation}
with the kernel
\[
K_H(t,s) = (t-s)^{\frac12-H} \FF\Bigl(\frac12-H,\; \frac12 - H,\;
\frac32-H,\; \frac{s-t}{t}\Bigr),
\]
where  $\FF(a,b,c,x)$ is the Gauss
hypergeometric function, and the constant $C_H$ is defined by
\[ 
C_H =
\biggl(\frac{\Gamma(2-2H)}{2H\Gamma(\frac12+H)(\Gamma(\frac32-H))^3}\biggr)^{\frac12},
\]
where $\Gamma(x)$ denotes the gamma function. See \cite{J06} for details of
the definition of the integral with respect to $B^H$ in formula \eqref{B
via BH}. Introduce the process $X_t$, $t\ge 0$, by
\[
X_t = C_H \int_0^t K_H(t,s) d Z_s.
\]
Using the above connection between $B_t$ and $B_t^H$ and computing the
corresponding integral with respect to $ds$ we obtain that
\[
X_t =  \theta L_H \int_0^t s^{\frac12-H}ds + B_t,
\]
and the filtrations of the processes $Z_t$ and $X_t$ coincide. The constant
$L_H$ is defined by
\[ 
L_H = \biggl(2H \Bigl(\frac32-H\Bigr)
B\Bigl(\frac12+H,\;2-2H\Bigr)\biggr)^{-\frac12},
\]
where $B(x,y)$ is the beta
function (see the appendix for computational details).

From the Cameron--Martin theorem or the Girsanov formula and the general
Bayes theorem (see e.g. Chapters~6,~7 in \cite{LS}, and the appendix) one can find that the
conditional distribution of $\theta$ is normal:
\begin{equation}
\Law(\theta \mid \F_t^X) =
N\biggl(\frac{a_t}{b_t},\; \frac{1}{b_t}\biggr)
\label{cond dist}
\end{equation}
with the processes
\[
a_t = \frac{\mu}{\sigma^2} + L_H\int_0^t s^{\frac12 -H} d X_s,\qquad b_t =
\frac{1}{\sigma^2} + \frac{L_H^2}{2-2H} t^{2-2H}.
\]
Then $R(\tau)$ can be
written as
\begin{equation}
R(\tau) = \E (\tau + h(a_\tau, b_\tau))
\label{R simple 1}
\end{equation}
with the function $h(a,b) = \min(\E \xi^+, \E\xi^-)$ for a normal random
variable $\xi\sim N(\frac ab, \frac1b)$. In the explicit form,
$h(a,b) = \frac{1}{\sqrt b} \phi(\frac{a}{\sqrt b}) -
\frac{|a|}{b}\Phi(-\frac{|a|}{\sqrt b})$,
where $\Phi,\phi$ are the standard normal distribution and density
functions.

Observe also that \eqref{cond dist} implies that for the optimal stopping
time $\tau^*$, the corresponding optimal function $d^*$ is equal to 1 if
$a_{\tau^*}>0$ and $-1$ if $a_{\tau^*}\le0$.

Define the process $\tilde B_t$, $t\ge 0$,
\[
\tilde B_t = X_t - L_H \int_0^t \E(\theta\mid \F_s^X) s^{\frac12 - H} ds = X_t
- L_H \int_0^t\frac{a_s}{b_s}  s^{\frac12 - H} ds.
\]
The innovation representation (see Chapter~7.4 in \cite{LS}) implies that
$\tilde B_t$ is a standard Brownian motion. Then the process $a_t$ satisfies
the SDE
\[
d a_t = \frac{a_t}{b_t} d b_t + L_H t^{\frac12-H}d\tilde B_t.
\]
Next we'll apply the It\^o formula to $h(a_\tau,b_\tau)$. In order to avoid
problems caused by that $h(a,b)$ is not smooth at $a=0$, consider the
function $\tilde h(a,b) = h(a,b) + \frac{|a|}{2b}$. It can be easily
verified that
$\tilde h \in C^{2,1}(\mathbb{R}\setminus\{0\}\times \mathbb{R}_+)$,  the derivative $\tilde h'_a(a,b)$
is continuous at $a=0$ for any $b> 0$, and 
$\tilde h'_b(a,b) + \frac12 \tilde h''_{aa}(a,b) + \frac {|a|}b \tilde
h'_a(a,b) = 0$
for all $a\in\R\setminus\{0\}$, $b> 0$. From this identity, applying the It\^o formula
we obtain
\[
d\tilde h (a_t,b_t) = \tilde h '_a(a_t, b_t) L_H t^{\frac12-H}d \tilde B_t.
\]
Using that
$\tilde h'_a(a,b) = - \frac{1}{b} \Phi(-\frac{|a|}{\sqrt b}) +
\frac{1}{2b}$,
one can see that
$\E \int_0^\infty (\tilde h'_a(a_t, b_t) t^{\frac12-H} )^{2} dt$ is finite,
so $\tilde h(a_t,b_t)$ is a square-integrable martingale.
Therefore, for any stopping time $\tau$ we have
$\E \tilde h(a_\tau,b_\tau) = \tilde h(a_0,b_0)$, which transforms \eqref{R
simple 1} to
\[
R(\tau) = \E\biggl(\tau - \frac{|a_\tau|}{2b_\tau}\biggr) + \tilde h(a_0,b_0).
\]
(Explicitly, $\tilde h(a_0,b_0) = \sigma \phi(\frac\mu\sigma) +
|\mu|(\frac12 - \Phi(-\frac{|\mu|}{\sigma}))$.)
Observe that $\frac{a_t}{b_t}$ satisfies the equation
\[
d\Bigl(\frac{a_t}{b_t}\Bigr) = L_H\frac{t^{\frac12-H}}{b_t}  d \tilde B_t.
\]
For brevity, denote
\[
\gamma = \gamma(H) = \frac{1}{2-2H}.
\]
Then under the following monotone change of time (see the appendix for details)
\begin{equation}
t(r) = \biggl(\frac{(2-2H)r}{\sigma^2 L_H^2(1-r)}\biggr)^\gamma, \qquad
r\in[0,1),
\label{time change}
\end{equation}
where $t$ runs through the half-interval $[0,\infty)$ when $r$ runs through
$[0,1)$, the process
\[
W_{r} = \frac{a_{t(r)}}{\sigma b_{t(r)}} - \frac{\mu}{\sigma}
\]
is a sBm in $r\in[0,1)$, and the filtrations
$\F_{r}^W$ and $\F_{t(r)}^X$ coincide. Therefore, for any stopping time $\tau$ of
the filtration $\F_t^X$ we have
\[
R(\tau) = \frac{\sigma}{2} \E\biggl(M_{\sigma,H}
\biggl(\frac{\rho}{1-\rho}\biggr)^\gamma - \Bigl|W_\rho +
\frac{\mu}{\sigma}\Bigr| \biggr) + \tilde h(a_0,b_0),
\]
where $\rho = t^{-1}(\tau)<1$ is a stopping time of the filtration
$\F_{r}^W$, and the constant
\[
M_{\sigma,H} = \frac{2}{\sigma}
\biggl(\frac{2-2H}{\sigma^2 L_H^2}\biggr)^\gamma.
\]
Thus, the optimal stopping problem for $X$ in $t$-time
is equivalent to the following optimal stopping problem for $W$ in $r$-time:
\begin{equation}
V = \inf_{\rho < 1} \E\biggl(M_{\sigma,H}
\biggl(\frac{\rho}{1-\rho}\biggr)^\gamma - \Bigl|W_\rho +
\frac{\mu}{\sigma}\Bigr|\biggr).\label{V}
\end{equation}
Namely, if $\rho^*$ is an optimal stopping time in \eqref{V}, then an optimal
decision rule $\delta^*=(\tau^*,d^*)$  is given by
\begin{equation}
\tau^* = t(\rho^*),\qquad d^* = \I(a_{\tau^*} > 0) - \I(a_{\tau^*} \le 0).
\label{trasnform}
\end{equation}

\section{The main result}
\label{main}
In this section we formulate the main theorem about the solution of problem
\eqref{V}, which provides an optimal sequential test via
\eqref{trasnform}. Throughout we will assume that the parameters of the
problem $\mu$, $\sigma$, $H$ remain fixed
and will denote the function
\[
f(t)  = M_{\sigma,H}\left(\frac{t}{1-t}\right)^\gamma.
\]

It is well-known that under general conditions the solution of an optimal
stopping problem for a Markov process can be represented as the first time
when the process enters some set -- a stopping set. Namely, let us first
rewrite our problem in the Markov setting by allowing the process $W_t$ to
start from any point $(t,x)\in[0,1)\times\mathbb{R}$:
\begin{equation}
V(t,x) = \inf_{\rho < 1-t} \E(f(t+\rho) - |W_\rho + x|) - f(t),\label{Vmarkov}
\end{equation}
where the infimum is over all stopping times $\rho$ of the Brownian motion
$W$ such that $\rho<1-t$ a.s. 
In particular, for the quantity $V$ from \eqref{V} we have
$V=V(0,\frac{\mu}{\sigma})$. We subtract $f(t)$ in the definition of
$V(t,x)$ to make the function $V(t,x)$ bounded. For $t=1$ we define
$V(1,x) = -|x|$.

The following theorem describes the structure of the optimal stopping time
in problem \eqref{Vmarkov}. In its statement, we set
\[
t_0 = t_0(H) := \max\biggl(0,\; 
\frac{1-2H}{4(1-H)}\biggr).
\]
Obviously, $0<t_0 < \frac14$ for $H<\frac12$ and $t_0=0$ for $H\ge \frac12$.
We find the solution of problem \eqref{Vmarkov} by describing the boundary
of the stopping set as a function of time $t$: in the case $H<\frac12$ this
will be done only for $t\ge t_0$, while in the case $H\ge \frac12$ for all
$t\in(0,1)$. Unfortunately, our method of proof does not work for $t<t_0$
when $H\le\frac12$, although the boundary can still be formally found from
the equation we obtain in the theorem; see the discussion below.

\begin{theorem}
1) There exists a function $A(t)$ defined on $(t_0,1]$, which is continuous,
non-increasing, and strictly positive for $t<1$ with $A(1)=0$, such
that for any $t> t_0$ and $x\in\mathbb{R}$ the optimal stopping time in 
problem \eqref{Vmarkov} is given by
\[
\rho^*(t,x) = \inf\{s\ge 0: |W_s+x| \ge A(t+s)\}.
\]
Moreover, for any $t\in(2 t_0,1]$ the function $A(t)$ satisfies the inequality
\begin{equation}
A(t) \le \frac{(1-t)^\gamma}{2M_{\sigma,H} t^{\gamma-1}}.
\label{A bound}
\end{equation}
2) The function $A(t)$ is the unique continuous non-negative solution of the
integral equation
\begin{equation}
\label{inteq}
G(t,A(t)) = \int_t^1 F(t,A(t),s,A(s)) ds, \qquad t\in(t_0,1),
\end{equation}
with the functions $G(t,x) = \E|\zeta\sqrt{1-t} + x| - x$ and $F(t,x,s,y) =
f'(s)\P(|\zeta\sqrt{s-t}+x| \le y)$ for a standard normal random variable $\zeta$. 
\end{theorem}

The main reason why we characterize the boundary $A(t)$ only for $t> t_0$ in
the case $H<\frac12$ is that the method of proof we use to show that $A(t)$
satisfies the integral equation requires it to be of bounded variation (at
least, locally). This condition is needed as a sufficient condition to apply
the It\^o formula with local time on curves \cite{P05}, on which the proof
is based on. In the case $H\ge\frac12$ and for $t\ge t_0$ in the case
$H<\frac12$ by a direct probabilistic argument we can prove that $A(t)$ is
monotone and therefore has bounded variation; this argument however doesn't
work for $t<t_0$ in the case $H<\frac12$, and, as a formal numerical
solution shows, the boundary $A(t)$ seems to be indeed not monotone in that
case. Of course, the assumption of bounded variation can be relaxed while
the It\^o formula can still be applied (see e.g. \cite{P05,E06,FPS}),
however verification of weaker sufficient conditions is problematic.
Although the general scheme to obtain integral equations of type
\eqref{inteq} and prove uniqueness of their solutions was discovered quite a
while ago (the first full result was obtained by Peskir \cite{P05b} for the
optimal stopping problem for American options), and has been used many times
in the literature for various optimal stopping problems (a large number of
examples can be found in \cite{PS}), we are unaware of any of its
non-trivial applications in the case when stopping boundaries are not
monotone. Nevertheless, a formal numerical solution of the integral equation
shows that the stopping boundaries ``look smooth'', but we admit that a
rigor proof of this fact in the above-mentioned cases remains an open
question.

Note also, that in the case $H\ge \frac 12$ the space-time transformation we
apply to pass from the optimal stopping problem for the process $a_t$ to
the problem for $W_r$ is essential from this point of view, because the
boundaries in the problem for $a_t$ are not monotone. Moreover, they are
not monotone even in the case $H=\frac12$, when $a_t$ is obtained by simply
shifting $X_t$ in time and space, see \cite{C1,ZM12}.

The second remark we would like to make is that in the case $H>\frac12$ we do
not know whether $A(0)$ is finite. In the case $H=\frac12$ the finiteness of
$A(0)$ follows from inequality \eqref{A bound}, which is proved by a direct
argument based on comparison with a simpler optimal stopping problem (one
can easily see from the proof that \eqref{A bound} extends to $t=0$ for
$H=\frac12$). It seems that a deeper analysis is required for the case
$H>\frac12$, which is beyond this paper.

\begin{figure}
\centering
\includegraphics[width=7cm]{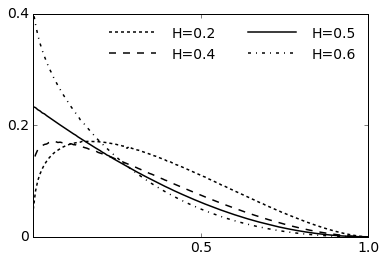}
\caption{The stopping boundary $A(t)$ for different values of $H$ and
$\sigma=1$.}
\label{fig1}
\end{figure}

Figure~\ref{fig1} shows the stopping boundary $A(t)$ for different values
$H$ computed by solving equation \eqref{inteq} numerically. The solution can
be obtained by backward induction on a discrete set of points of $[t_0,1]$
starting with $t=1$ and going towards $t=t_0$ using that the expression
under the integral depends only on the values of $A(s)$ for $s\ge t$; the
method is described in more details, for example, in \cite{PP}.

\section{The proof of the theorem}
\label{proofs}
The proof will be conducted in several steps: (i) prove the monotonicity and
continuity of the function $V(t,x) +|x|$; (ii) analyze the stopping set and
its boundaries; (iii)
formulate a free-boundary problem for $V(t,x)$; (iv) derive the integral
equation for $A(t)$; (v) prove the uniqueness of its solution.

\medskip
(i) It is clear that
\[
V(t,x) = \inf_{\rho < 1-t} \E\biggl(\int_0^\rho f'(s+t)ds - |W_\rho +
x|\biggr).
\]
One can easily verify that $f'(t)$ is increasing in $t\in[t_0,1)$, which
implies that $V(t,x)$ doesn't decrease in $t\in[t_0,1]$ for each fixed $x$
($t_0$ is the point of minimum of $f'(t)$ on $(0,1)$).
Also, for any $t\in[0,1]$ the function $V(t,x)+|x|$ is non-decreasing in $x$
on $[0,\infty)$ and non-increasing in $(-\infty,0]$. This follows from the
inequality $|x_1| - \E|W_\rho+x_1| \le |x_2| - \E|W_\rho+x_2|$ for any
$0\le x_1 \le x_2$ and any stopping time $\rho<1-t$.

The monotonicity of $V(t,x)+|x|$ both in $t$ and $x$ implies that in order
to prove the continuity of this function for $t\ge t_0$, it is enough to
prove the continuity in each argument (see e.g. \cite{KD69}). The continuity
in $x$ follows from that for any $t\in[0,1)$, $x_1,x_2\in\mathbb{R}$ we have
\begin{equation}
\bigl|V(t,x_1) - V(t,x_2)\bigr| \le \sup_{\rho< 1-t} \E\bigl| |W_\rho + x_1| -
|W_\rho+x_2|\bigl| \le |x_1-x_2|,
\label{dV bound}
\end{equation}
and also the continuity of $V(1,x)$ is obvious. In order to prove the
continuity in $t$, fix $x\in\mathbb{R}$ and consider arbitrary $t_0\le t<1$,
$0\le\epsilon<1-t$. Then for any stopping time $\rho<1-t$
\begin{equation}
\begin{split}
&\E\left[ \int_0^\rho f'(t+s) - |W_\rho+x| \right]  \\ &\qquad\ge
\E\biggl[ \int_0^{\rho\vee\epsilon} f'(t+s) ds - |W_{\rho\vee\epsilon} + x|
+|W_\epsilon+x|-|W_{\rho\wedge \epsilon} +x| \biggr] - \int_0^\epsilon f'(t+s)ds \\
&\qquad\ge \E V(t+\epsilon, W_\epsilon+x) +
\E(|W_\epsilon+x|-|W_{\rho\wedge\epsilon}+x|),
\end{split}
\label{2.c.eq:2}
\end{equation}
where the first inequality follows from straightforward algebraic
transformations, and to prove the second one the Markov property of Brownian
motion was used. We have
\[
\bigl|\E(|W_\epsilon+x| - |W_{\rho\wedge \epsilon}+x|)\bigr| \le
\E|W_{\epsilon} - W_{\rho\wedge\epsilon}| = \E|\widetilde
W_{\epsilon-\rho\wedge\epsilon} |\le \sqrt{\epsilon}, 
\]
where the process $\widetilde W=(\widetilde W_t)_{t\le \epsilon}$,
$\widetilde W_t = W_\epsilon- W_{\epsilon-t}$, is a Brownian motion. In the
right inequality we used the well-known bound
$\E |\widetilde W_\tau| \le \sqrt{\E \tau}$, which easily follows from the
Wald identity and Jensen's inequality. Taking in \eqref{2.c.eq:2} the
infimum over all stopping times $\rho< 1-t$ we obtain
$V(t,x) \ge \E V(t+\epsilon, W_\epsilon+x) -\sqrt{\epsilon}$. According to
what was proved above,
$|\E V(t+\epsilon,W_\epsilon+x) - V(t+\epsilon,x)| \le \E|W_\epsilon|
=\sqrt{2\epsilon/\pi}$.
Therefore,
$0\ge V(t,x) - V(t+\epsilon,x) \ge -(\sqrt{\epsilon} +
\sqrt{2\epsilon/\pi})$,
which proves the continuity of $V(t,x)$ in $t$ for $t\in[t_0,1)$. The
continuity at $t=1$ is obvious since
$V(1,x) = -|x| \ge V(t,x) \ge -\sup_{\rho\le 1-t}\E|W_\rho+x| \to -|x|$ as
$t\to1$.

\medskip (ii) Define the stopping set $D=\{(t,x) \in [0,1]\times\mathbb{R}
\mid V(t,x) = -|x|\}$ and the set $D_0 = D\cap \{t\ge t_0\}$. The continuity
of $V(t,x)$ implies that $D_0$ is closed.

It is clear that $D_0$ is symmetric in $x$, and from the monotonicity of
$V(t,x)+|x|$ it follows that $D_0$ can be represented in the form
\[
D_0 = \{(t,x) : t\in[t_0,1],\; |x| \ge A(t)\},
\]
where $A(t)$ is some non-increasing function on $[t_0,1]$. Obviously,
$A(1)=0$ and one can easily see that $A(t)>0$ for any $t<1$, since for for
any $t<1$ it is possible to find a sufficiently small non-random $r$ such
that $\E(f(r+t) - |W_r|) < f(t)$, hence $V(t,0)< 0$.

Let us prove inequality \eqref{A bound} for $A(t)$. By using the inequality
$f'(s) \ge  \frac{M_{\sigma,H} \gamma t^{\gamma-1}}{(1-s)^{\gamma+1}}$ for any
$s\ge t$, one can see that for any
$t\in[0,1)$ and any stopping time $\rho <1-t$
\[
f(t+\rho) - f(t) \ge \frac{M_{\sigma,H} t^{\gamma-1}}{(1-t)^{\gamma}} 
\left( \frac{(1-t)^\gamma}{(1-\rho-t)^\gamma} -1\right).
\]
Denote the expression in the brackets by $\nu(\rho)$. Observe that if $t\ge
2t_0 = 1-\gamma$, then $\nu(\rho) \ge
\rho$  and hence $\nu$ is a stopping time, and also
$\E|W_\nu+x| \ge \E|W_\rho+x|$. Then
\[
V(t,x) \ge \inf_{\nu \ge 0} \E
(c\nu - | W_\nu + x|),
\]
where $c = \frac{M_{\sigma,H} t^{\gamma-1}}{(1-t)^{\gamma}}$. For the
optimal stopping problem in the RHS, the solution is well-known (see e.g.
Section~16 in~\cite{PS}): the optimal stopping time is
$\nu^*(x)=\inf\{t\ge 0 : |W_t+x| \ge \frac1{2c}\}$. Hence, if
$|x|\ge\frac{1}{2c}$, then $V(t,x) = -|x|$ and $(t,x)\in D$, which proves
\eqref{A bound}. In particular, $A(t)$ is finite-valued for all $t>2t_0$.

In order to show that $A(t)$ is finite-valued for $t\in[t_0,2t_0]$ when $H<\frac
12$, one can use that $V(t,x) \ge \inf_{\rho}\E (f'(t_0) \rho -
|W_\rho + x|)$, since $t_0$ is the point of minimum of the function $f'(t)$
when $H<\frac 12$. Again, all the points $(t,x)$ with $|x|\ge
\frac{1}{2f'(t_0)}$ should be in the stopping set, so $A(t)$ is bounded by $1/(2f'(t_0))$.

Next, for any point $(t,x)$, $t> t_0$, define the candidate optimal
stopping time $\rho^*(t,x)$ -- the first entry into the stopping set:
\[
\rho^* = \rho^*(t,x) =  \inf\{s\ge 0 : (t+s,W_s+x) \in D\}.
\]
In the general theory of optimal stopping for Markov processes, it is
well-known that the first entry into the stopping set is an optimal stopping
time under mild conditions. In our problem this fact can be proved
similarly to \cite{ZS} (one subtlety here is that general conditions for the optimality of
$\rho^*$ typically require some boundedness of the payoff function, see for
example Chapter 1 in \cite{PS}; since the payoff in our problem is
unbounded a finer argument is needed).


Now we prove that $A(t)$ is continuous on $(t_0,1]$. Since it doesn't
increase and the set $D$ is closed, it is clear that $A(t)$ is
right-continuous. Let us prove that it is left-continuous. Using that
$\E W_{\rho^*} = 0$ we can write
\begin{equation}\label{eq:8}
V(t,x) = \E \left( f(\rho^*+t) + W_{\rho^*} - |W_{\rho^*} + x| 
\right) - f(t).
\end{equation}
Suppose $A(t-)> A(t)$ for some $t\in(t_0,1)$. Consider points $(t-\varepsilon, x)$ with $x =
(A(t-)+A(t))/2$ and sufficiently small $\varepsilon>0$. Let
$\Omega_\varepsilon=\{\omega:W_s(\omega)+x \in [A(t),A(t-)]\text{ for all } s\le
\varepsilon\}$ denote the random event that $W$ exists
the rectangle $[t-\epsilon,t]\times[A(t), A(t-)]$ through the right boundary.
Denote for brevity $f_\epsilon = f(t)-f(t-\epsilon)>0$. Then from
\eqref{eq:8} we have
\[
\begin{split}
V(t-\varepsilon,x) &\ge \E\bigl[ (f_\epsilon - x)
\I(\Omega_\varepsilon)\bigr] + \E\bigl[\bigl(W_{\rho^*} - |W_{\rho^*} +
x|\bigr)\I(\Omega\setminus A_\varepsilon)\bigr] \\ &\ge
f_\epsilon-x -  \P(\Omega\setminus A_\varepsilon) \left(f_\varepsilon-x +
\sqrt{\E(W_{\rho^*} - |W_{\rho^*}+x|)^2} \right),
\end{split}
\]
where in the first inequality we used that
$W_{\rho^*} - |W_{\rho^*}+x| = -x$ on $\Omega_\epsilon$ and in the second
one we applied the Cauchy--Schwarz inequality to
$\E[(W_{\rho^*} - |W_{\rho^*} + x|)\I(\Omega\setminus A_\varepsilon)]$.

According to Doob's martingale inequality,
$\P(\Omega\setminus \Omega_\varepsilon) \le
2\exp\left(-\frac{(A(t-)-A(t))^2}{8\varepsilon}\right)$. Since
$f_\varepsilon\ge f'(t-\epsilon)\epsilon$ when $t-\epsilon\ge t_0$, there
exists a sufficiently small $\varepsilon>0$ such that
$V(t-\varepsilon,x) > -x$, which contradicts the definition of $V(t,x)$.
This proves the continuity of $A(t)$ on $(t_0,1)$. The continuity at $t=1$
follows from inequality \eqref{A bound}.

\medskip (iii) As follows from the general theory of optimal stopping for
Markov processes, inside the continuation set
$C_0= (t_0,1]\times \R \setminus D$ the value function $V(t,x)$ is $\C^{1,2}$
and satisfies the partial differential equation (see Section 7 in \cite{PS})
\begin{equation}
V'_t(t,x) + \frac12 V''_{xx}(t,x) = - f'(t), \qquad t\in(t_0,1),\; |x| <
A(t).
\label{V pde}
\end{equation}
Together with the condition $V(t,x) = -|x|$ in the set $D$, this constitutes
a free boundary problem for the value function $V(t,x)$ with the unknown
free boundary $A(t)$. 

The continuity of $V(t,x)$ implies the so-called condition of continuous fit: $V(t,A(t)-) =
V(t,A(t)+)$, i.e. $V(t,x)$ is continuous at the stopping boundary. Let us now prove
the smooth-fit condition, which states that the $x$-derivative of $V(t,x)$
is continuous at the stopping boundary:
\begin{equation}
V'_x(t,A(t)-) = V'_x(t,A(t)+) \; (=-1) , \qquad t\in(t_0,1].
\label{smooth fit}
\end{equation}
The function $V(t,x)$ in concave in $x$ since it is the infimum (over
$\rho$) of concave functions. Therefore, there exist the left and right
derivatives $V'_{x\pm}(t,A(t))$. Clearly, $V'_{x+}(t,A(t)) = -1$, since
$V(t,x) = -|x|$ for $x\ge A(t)$. Moreover, for any sufficiently small
$\epsilon>0$ we have
\[
\frac{V(t,A(t)-\epsilon) - V(t,A(t))}{-\epsilon} \ge -1
\]
since $V(t,A(t)-\epsilon) \le -(A(t)-\epsilon)$ and $V(t,A(t)) = -A(t)$.
Therefore, $V'_{x-}(t,A(t)) \ge -1$. Let us prove the opposite inequality.

Fix $t\in(t_0,1)$. Set $x=A(t)$ and let $\epsilon>0$ be sufficiently small.
Then for the optimal time $\rho^* = \rho^*(t,x-\epsilon)$ we have
\[
\begin{aligned}
\frac{V(t,x-\epsilon)- V(t,x)}{-\epsilon} \le \frac{
\E |W_{\rho^*} + x| - \E|W_{\rho^*} + x-\epsilon|}{-\epsilon},
\end{aligned}
\]
where we used that $V(t,x) \le \E\bigl[ f(\rho^*+t) - |W_{\rho^*}+x| \bigr]
- f(t)$ and $V(t,x-\epsilon) = \E\bigl[ f(\rho^*+t) - |W_{\rho^*}+x-\epsilon| \bigr] - f(t)$.
Transform the obtained expression:
\begin{multline*}
\E |W_{\rho^*} + x| - \E|W_{\rho^*} + x-\epsilon|  = \epsilon
\P(W_{\rho^*}+x-\epsilon = A(t+\rho^*)) \\ +
\E\bigl[\bigl(|W_{\rho^*}+x|-|W_{\rho^*}+x-\epsilon|\bigl) \I(W_{\rho^*}+x-\epsilon = -A(t+\rho^*) \bigr].
\end{multline*}
The second term can be bounded in absolute value by $\epsilon
\P(W_{\rho^*}+x-\epsilon = -A(t+\rho^*)) = o(\epsilon)$. Then
\[
\lim_{\epsilon\downarrow 0} \frac{V(t,x-\epsilon)- V(t,x)}{-\epsilon} \le -1,
\]
which proves the inequality $V'_{x-}(t,A(t))\le -1$.

\medskip
(iv) So far we have established the following properties: (a)~$V(t,x)$ is continuous on
$[t_0,1]\times\R$ and is
$\C^{1,2}$ in $C_0$ and in the interior of $D_0$; (b)~$A(t)$ is continuous and non-increasing on
$(t_0,1]$; (c)~$( V'_t + \frac12 V''_{xx})(t,x)$ is locally bounded in~$C_0$ and
in the interior of~$D_0$, which follows from \eqref{V pde};
(d) the function~$x \mapsto V(t,x)$ is concave and
the function $t\mapsto V'_x(t,A(t)\pm)\;(\equiv\mp1)$ is continuous.

These properties allow to apply the
It\^{o} formula with local time on curves~(see \cite{P05} and
Section~2.6 in~\cite{P05b}) to~$V(t,x)$:
for any $t_0<t<T<1$ and $x\in\R$ we have
\begin{equation}
\label{V Ito}
\begin{split}
&\E V(T,W_{T-t}+x) - V(t,x) \\ &\quad= \E \int_0^{T-t} \Bigl(V'_t + \frac12 V''_{xx}\Bigr)(t+s,W_s+x) \I(W_s+x \neq \pm A(t+s))\,ds\\ 
&\quad+\,\E \int_0^{T-t} V'_x(t+s,W_s+x) \,\I(W_s+x \neq \pm A(t+s))\,d
W_s\\ 
&\quad+\,\frac12 \E \int_0^{T-t} \Delta V'_x(t+s,A(t+s)) \,\I(W_s+x=A(t+s))\,d
L^{A}_s\\ 
&\quad+\,\frac12 \E \int_0^{T-t} \Delta V'_x(t+s,-A(t+s)) \,\I(W_s+x=-A(t+s))\,d
L^{-A}_s,
\end{split}
\end{equation}
where $L^{\pm A}$ is the local time processes of~$W$ on the curves~$\pm A$
(see~\cite{P05}), and $\Delta V_x'(t,x) = V_x'(t,x+) - V_x'(t,x-)$.

Smooth-fit condition~\eqref{smooth fit} implies that
the two last terms in~\eqref{V Ito} are equal to zero. Also,
the derivative $V'_x(t,x)$ is uniformly bounded according to \eqref{dV
bound}, and therefore the expectation of the stochastic integral
in~\eqref{V Ito}
is also zero.

From \eqref{V pde} and the fact that $(V'_t+\frac12 V'_{xx})(t,x) = 0$ for
$|x|>A(t)$  we obtain
\begin{equation}
\label{V repr}
V(t,x) = \E V(T,W_{T-t}+x) + \int_0^{T-t}  f'(t+s)\P(|W_s+x|
< A(t+s))\,ds.
\end{equation}
By passing to the limit $T\to 1$, we have
$\E V(T,W_{T-t}+x) \to -\E |W_{1-t}+x|$ from the dominated convergence
theorem. Finally, to obtain integral equation \eqref{inteq}, it remains to
put $x= A(t)$ and use the identity $V(t,A(t)) = -A(t)$.

(v) To prove that $A(t)$ is the  unique solution of  integral equation
\eqref{inteq}, suppose $\tilde A(t)$ is another non-negative continuous
solution satisfying \eqref{A bound}. Define
the function (cf. \eqref{V repr})
\[
\begin{split}
\tilde V(t,x) &= -\E |W_{1-t}+x| + \E \int_0^{1-t} f'(t+s) \I(|W_{s} + x| <
\tilde A(t+s)) ds \\
&= -|x|-G(t,x) +  \int_t^1 F(t,x,s,\tilde A(s))ds,
\end{split}
\]
where $ t\in(t_0,1)$, $x\in \R$. We don't exclude the possibility $\tilde
V(t,x)=+\infty$, however, obviously, $\tilde V(t,x)>-\infty$ for all $t,x$. Using the strong Markov property, one can
show that for any  $t\in(t_0,1)$, $x\in\R$, and any stopping time  $\rho<
1-t$ we have
\begin{equation}\label{1.13}
\begin{aligned}
\tilde V(t,x) = \E \tilde V(t+\rho, W_\rho+x) + \E  \int_0^\rho
f'(t+s) \I(|W_{s} + x|
< \tilde A(t+s)) ds.
\end{aligned}
\end{equation}
Consider the stopping time
$\rho_{\tilde A} =\inf\{s\ge 0: |W_s+x| = \tilde A(t+s)\}\wedge (1-t)$.
Since $\tilde A$ satisfies the integral equation, one can see that
$\tilde V(t+\rho_{\tilde A}, W_{\rho_{\tilde A}} +x)= -|W_{\rho_{\tilde
A}}+x|$. Together with \eqref{1.13}, this implies
\begin{align*}
&\tilde V(t,x) = -|x|,& \quad &|x|\ge \tilde A(t),\\
&\tilde V(t,x) = \E\biggl[  \int_0^{\rho_{\tilde A}} f'(t+s)  ds -
  |W_{\rho_{\tilde A}}+x| \biggr],&
&|x| < \tilde A(t).
\end{align*}
Consequently $\tilde V(t,x)\ge V(t,x)$ for all $t\in(t_0,1]$, $x\in \R$.

Suppose $\tilde A(t) > A(t)$ for some $t\in(t_0,1)$. Set $x= \tilde A(t)$
and consider the corresponding optimal stopping time $\rho^*$. Then from
\eqref{1.13}, using that $\tilde V(t,x)=-|x|$ and
$\tilde V(t+\rho^*,W_{\rho^*})\ge V(t+\rho^*,W_{\rho^*})= -|W_{\rho^*}+x|$,
we get
\[
-x \ge 
\E\biggl[-|W_{\rho^*}+x| + \int_0^{\rho^*}
f'(t+s) \I(|W_s+x| < \tilde A(t+s))ds\biggr].
\]
However, the expectation of the integral in the above formula is strictly positive since the
process $W_s+x$ spends a.s. strictly positive time between the boundaries
$\pm \tilde A$ (here we use the assumption that $\tilde A$ is continuous). Moreover, $\E |W_{\rho^*} +x| = \E(W_{\rho^*}+x) = x$, since
$W_s+x$ remains positive until time $\rho^*$. Thus, we get a
contradiction, implying that $\tilde A(t) \le A(t)$.

Suppose now that $\tilde A(t) < A(t)$ for some
$t\in(t_0,1)$ and set  $x =
\tilde A(t)$. Then $\tilde V(t,x) = -x$ and $\tilde V(t+\rho^*, W_{\rho^*}) = -
|W_{\rho^*} + x|$. From  \eqref{1.13} we get
\[
-x = 
\E\biggl[-|W_{\rho^*}+x| + \int_0^{\rho^*}
f'(t+s) \I(|W_s+x| < \tilde A(t+s))   ds\biggr]
< V(t,x),
\]
where we used that the indicator function under the integral is not
identically 1 a.s. Again, we get a contradiction with that $V(t,x) \le
-x$, which finishes the proof.

\section*{Appendix}
This appendix contains some technical details related to Section~\ref{reduction}.

\medskip
(a) Evaluation of the integral $C_H\int_0^t K_H(t,s)ds$. By the change $u=\frac{t-s}{s}$
we have
\begin{equation}
\int_0^t K_H(t,s)ds = t^{\frac32-H}\int_0^\infty \frac{u^{\frac12 -
H}}{(u+1)^{\frac52 -H}} \;\FF\Bigl(\frac12-H,\; \frac12-H,\; \frac32-H,\;
-u\Bigr)du.
\label{A.1}
\end{equation}
Next one can use the formula (see \cite{PB}, formula 2.21.1.16)
\begin{multline*}
\int_0^\infty \frac{u^{c-1}}{(u+z)^{r}}\; \FF(a,b,c,-u)dx \\=
\frac{\Gamma(c)\Gamma(a-c+r)\Gamma(b-c+r)}{\Gamma(r)\Gamma(a+b-c+r)}
\FF(a-c+r,\;b-c+r,\;a+b-c+r,\;1-z),
\end{multline*}
which holds for complex $a,b,c,r,z$ such that $\Re(a+r),\Re(b+r)>\Re c > 0$
and $|\arg z| < \pi$. Clearly, these conditions are satisfied in \eqref{A.1},
and using that $\FF(a,b,c,0) = 1$ we obtain
\[
\int_0^t K_H(t,s) ds =  \frac{t^{\frac32 -H} \Gamma(\frac32 -H)^2}{(\frac32
-H) \Gamma(2-2H)}.
\]
Multiplying the right-hand side by the constant $C_H$, we obtain $L_H$.

\medskip
(b) The conditional distribution $\Law(\theta \mid \F_t^X)$. Although
formula \eqref{cond dist} can be found from the general filtration theory
for Gaussian processes, let us show its straightforward derivation.
From the general Bayes theorem (see e.g. \S\,II.7 in \cite{S}), one finds the conditional density
\[
\P(\theta \in du \mid \F_t^X) = {\dfrac{d \P_t^u}{d\P_t^0}
\phi_{\mu,\sigma}(u)} \times \biggl({\displaystyle \int_{\R}\dfrac{d \P_t^v}{d\P_t^0}
\phi_{\mu,\sigma}(v) dv}\biggr)^{-1}, \qquad u\in\R,
\]
where $\frac{d \P_t^u}{d\P_t^0}$ denotes the density process of the 
measure generated  by $X^u_t = B_t + uL_H \int_0^t
s^{\frac12-H}ds$ with respect to the measure generated by $X_t^0$, both
restricted to the $\sigma$-algebra $\F_t^X$. By $\phi_{\mu,\sigma}(u) =
\frac{1}{\sigma}\phi(\frac{u-\mu}{\sigma})$ we denote the density function
of the normal distribution $N(\mu,\sigma^2)$. From the Cameron--Martin
theorem
\[
\frac{d \P_t^u}{d\P_t^0} = \exp \biggl(uL_H  \int_0^t  s^{\frac12 -H} d X_s -
\frac{u^2L_H^2}2 \int_0^t s^{1-2H} ds \biggr)
\]
and the remaining step to obtain $\Law(\theta \mid \F_t^X)$ is a
straightforward integration.

\medskip (c) The change of time for $a_t/b_t$. The Dambis--Dubins--Schwarz
theorem (see \cite{RY}) says that any continuous local martingale $Y_t$ with
$Y_0=0$ can be obtained from a Brownian motion $W_r$ by the time change
$Y_t = W_{r(t)}$, where $r(t) = \langle Y\rangle_t$, the quadratic
variation of $Y$. This representation is strict in the sense that $W$ is
defined on the same probability space as $Y$. In our case
\[\Bigl\langle \frac{a}{\sigma b}\Bigr\rangle_t = \int_0^t 
\frac{L_H^2 s^{1-2H}}{\sigma^2 b_s^2} ds = 1 - \biggl(1+\frac{\sigma^2L_H^2
t^{2-2H}}{2-2H}\biggr)^{-1},
\]
which gives formula \eqref{time change} for $t(r)$.

\end{document}